\documentclass{article}

\usepackage{cmap}			 % поиск в PDF
\usepackage{mathtext} 		 % русские буквы в формулах
\usepackage[T2A]{fontenc}	 % кодировка
\usepackage[utf8]{inputenc}	 % кодировка исходного текста

\usepackage[inline]{asymptote}	
\usepackage{svg}
\usepackage{enumerate}
\usepackage{fancyhdr}
\usepackage{titlesec}
\usepackage[letterpaper,top=2cm,bottom=2cm,left=3cm,right=3cm,marginparwidth=1.75cm]{geometry}
\usepackage{amsmath,amsfonts,amssymb,amsthm,mathtools,tikz-cd,circuitikz,tikz}
\usepackage{leftindex}
\usepackage{graphicx}
\usepackage[colorlinks=true, allcolors=blue]{hyperref}
\usepackage{fancybox,fancyhdr}
\usepackage{xcolor}
\usepackage{faktor}
\definecolor{sangria}{rgb}{0.57, 0.0, 0.04}

\DeclareMathOperator{\Ker}{Ker}
\DeclareMathOperator{\im}{im}
\DeclareMathOperator{\Img}{Im}

\DeclareMathOperator{\coker}{coker}
\DeclareMathOperator{\Coker}{Coker}
\DeclareMathOperator{\coim}{coim}
\DeclareMathOperator{\Coim}{Coim}

\title{\textbf{Classes of Universal Epi- and Monomorphisms in Quasi-Abelian Categories}}
\author{Max Zinchenko}
\date{}
\begin{document}
	\maketitle
\begin{abstract}
	This paper provides a solution to the open problem formulated in \cite{D}, as well as examples of non-strict universal epimorphisms and monomorphisms.
\end{abstract}

\section*{Introduction}
The paper \cite{D} describes a theorem of fundamental importance in homological algebra, also known as the snake lemma for quasi-abelian categories, which requires the vertical arrows to be strict (see also \cite{H} for a non-additive setting). However, in \cite{F} this requirement is relaxed, demanding instead the so-called universality, a concept which will be discussed throughout the present work. Furthermore, in their work, Glotko and Kuzminov introduce the following axioms:
\par \textbf{Axiom 1:} Universal monomorphisms are stable under pullback.
\par \textbf{Axiom 1}$^{\ast}$\textbf{:} Universal epimorphisms are stable under pushout.
\par It was noted that whether these axioms hold in an arbitrary quasi-abelian category remains an open question. Although \cite{F} generalizes \cite{D}, it provides no example of a universal epi- or monomorphism in a non-integral category that is not strict, nor any examples of categories satisfying either of these axioms.
\section{Preabelian and quasi-abelian categories}
In this section, we will recall the basic definitions that we will need later. Recall that an additive category is called \textit{preabelian} if every morphism has a kernel and a cokernel. \\
Every morphism in preabelian category admit a canonical decomposition:
$$\begin{tikzcd}
	\Ker f \arrow[r, "\ker f", tail] & A \arrow[r, "f"] \arrow[d, "\coim f", two heads] & B \arrow[r, "\coker f", two heads] & \Coker f \\
	& \Coim f \arrow[r, "\overline{f}"]                & \Img f \arrow[u, "\im f", tail]     &         
\end{tikzcd}$$
morhism $f$ is a strict if $\overline{f}$ isomorphism.
Note that $f$ is a strict monomorphism (epimorphism) if and only if $f$ is a kernel (cokernel).\\
In the rest of this section, the category is assumed to be preabelian.
\par Let $\mathcal{A}$ be a category. We say that a class of morphisms $\mathcal{X}$ is \textit{stable under pushout} if, for any arbitrary morphism $g$ in $\mathcal{A}$ and any morphism $f \in \mathcal{X}$ in the pushout square below
$$\begin{tikzcd}
	X \arrow[r, "f"] \arrow[d, "g"'] \arrow[r, "\text{PO}", phantom, bend right=60] & Y \arrow[d] \\
	Z \arrow[r, "f'"']                                                              & P          
\end{tikzcd}$$
$f' \in \mathcal{X}$. Being \textit{stable under pullback} is defined dually.\\
\par If cokernels in $\mathcal{A}$ are stable under pullback, then $\mathcal{A}$ is called \textit{left quasi-abelian}. Dually, if kernels in $\mathcal{A}$ are stable under pushout, then $\mathcal{A}$ is called \textit{right quasi-abelian}. If $\mathcal{A}$ is both left and right quasi-abelian, then it is called \textit{quasi-abelian}.
\par If epimorphisms in $\mathcal{A}$ are stable under pullback, then $\mathcal{A}$ is called \textit{left integral}. Dually, if monomorphisms in $\mathcal{A}$ are stable under pushout, then $\mathcal{A}$ is called \textit{right integral}. If $\mathcal{A}$ is both left and right integral, then it is called \textit{integral}.
\par \textbf{Proposition 1.1}\cite[p. 173]{R} A quasi-abelian category is left integral if and only if it is right integral.

\section{Torsion theories}
Let $\mathcal A$ be a preabelian category. A \textit{torsion pair} in $\mathcal A$ is an ordered pair $(\mathcal T, \mathcal F)$ of full subcategories of $\mathcal A$ satysfying the following:
\begin{enumerate}
	\item $\text{Hom}_{\mathcal A}(\mathcal T, \mathcal F)=0$;
	\item for all $M$ in $\mathcal A$ there exists a short exact sequence
	$$\begin{tikzcd}
		0 \arrow[r] & A \arrow[r, "i_M", tail] & M \arrow[r, "p_M", two heads] & B \arrow[r] & 0
	\end{tikzcd}$$
	with $A \in \mathcal T$, $B \in \mathcal F$ and $i_M$, $p_M$ is strict mono and epimorphism respectively;
\end{enumerate}
By \cite[Proposition 2.5]{P}, there exists a functor $F : \mathcal{A} \to \mathcal{T}$ that is an idempotent and radical kernel subfunctor of the identity such that $\mathcal{T} = \{M \in \mathcal{A} \mid F(M) \cong M \}$. Moreover, $F$ is right adjoint to the canonical inclusion $\mathcal{T} \hookrightarrow \mathcal{A}$;
hence $F$ preserves limits.
\par It follows that pushouts and pullbacks in $\mathcal{T}$ can be computed via the ambient abelian
category $\mathcal{A}$. Indeed, let $f\colon A\to B$ and $g\colon A\to X$ be morphisms in $\mathcal T$,
and consider their pushout in $\mathcal A$:
$$\begin{tikzcd}
	A \ar[r," f"] \ar[d,"g"'] \ar[r, phantom, shift right=4ex, "{\rm PO}" marking] & B \ar[d,""] \\
	X \ar[r,"f'"] & Y
\end{tikzcd}$$
By the definition of a pullback, $Y:= \Coker(A \xrightarrow{[f, -g]} B \oplus X)$. By \cite[Proposition 2.9]{P}, $B \oplus X \in \mathcal{T}$ and therefore $Y \in \mathcal{T}$ as well. Now let $f\colon A\to B$ and $g\colon X\to B$ be morphisms in $\mathcal T$, and consider their
pullback in $\mathcal A$:
$$\begin{tikzcd}
	D \ar[r," f'"] \ar[d,""'] \ar[r, phantom, shift right=4ex, "{\rm PB}" marking] & X\ar[d,"g"] \\
	A \ar[r,"f"] & B
\end{tikzcd}$$
In general, $D$ need not belong to $\mathcal T$ since torsion classes are not closed under kernels.
However, because $F$ preserves limits, $F(D)$ is the pullback of $f$ and $g$ in $\mathcal T$.
$$\begin{tikzcd}
	F(D) \ar[r," Ff'"] \ar[d,""'] \ar[r, phantom, shift right=4ex, "{\rm PB}" marking] & X\ar[d,"g"] \\
	A \ar[r,"f"] & B
\end{tikzcd}$$
and by \cite[Proposition 2.5]{P}, we have $Ff' = f' \circ \mu$, where $\mu : F(D) \rightarrowtail D$.

\section{Definition and existence of universal epi- and monomorphisms}
\par We adopt the convention that any category is assumed to be quasi-abelian unless stated otherwise.
\par A monomorphism is called a \textit{universal monomorphism} if it is stable under pushout.
Dually, a \textit{universal epimorphism} is defined. The classes of a universal monomorphisms and universal epimorphisms are denoted by $M_u$ and $P_u$, respectively. 
\par It is obvious that $M \subseteq M_u \subseteq M_c$, where $M$ denotes the class of all monomorphisms and $M_c$ the class of strict monomorphisms. An analogous chain of inclusions holds for universal epimorphisms. Note that the equality $M = M_u$ holds if and only if the category is integral. Moreover, there exist examples of quasi-abelian categories where the proper inclusion $M \subsetneq M_u = M_c$ holds, and of course, in abelian categories all three inclusions is equalities.\\ 
We say that a morphism $f$ belongs to the class $MO_u$ if the morphism $\overline{f}$ from the canonical decomposition of $f$ belongs to the class $M_u$. Dually, the class $PO_u$ is defined.

\par \textbf{Proposition 3.1.} There exists a non-integral quasi-abelian category in which there exists a universal monomorphism that is not strict.\\
\textit{Proof.} Consider the category of normed spaces $\mathbf{Norm}$; its non-integrality and quasi-abelianity are described in \cite[Proposition 3.2.17]{B} and \cite[Theorem 3.2]{A}. Its objects are normed spaces, and its morphisms are continuous linear maps.
\par Let $A$ be an incomplete normed space and $A \xrightarrow{i} \hat{A}$ be its completion morphism. The morphism $i$ is a bimorphism, which is moreover not strict, as $i(A)$ is not closed in $\hat{A}$. Let $A \xrightarrow{g} G$ be an arbitrary morphism in $\mathbf{Norm}$, and consider the following pushout:
$$\begin{tikzcd}
	A \arrow[rr, "i", hook] \arrow[dd, "g"'] \arrow[rr, "\text{PO}", phantom, bend right=60] &                                 & \hat{A} \arrow[d, "{(0,id)}"]   \\
	&                                 & G \oplus \hat{A} \arrow[d, "p"] \\
	G \arrow[r, "{(id,0)}"]                                                                  & G \oplus \hat{A} \arrow[r, "p"] & D                              
\end{tikzcd}$$
Here $D := \faktor{G \oplus \hat{A}}{\overline{ran\begin{bsmallmatrix}
			-g \\
			i\\
\end{bsmallmatrix}}}$, $p$ is the canonical projection, and for a mapping $f \colon X \to Y$, $ran(f)=f(X)$.
\par Let us denote $j :=p \circ (id,0)$. By definition $j(m) = p((m,0))$, i.e $j(m)=0$ if $(m,0) \in \overline{ran\begin{bsmallmatrix}
		-g \\
		i\\
\end{bsmallmatrix}}$.
\par Let us prove that no such pairs with $m \neq 0$ exist there. Suppose there exists a sequence $\{a_n\}$, such that $i(a_n) \xlongrightarrow{n \rightarrow \infty} 0$ and $-g(a_n) \xlongrightarrow{n \rightarrow \infty} m$. Since $i(a_n) = a_n$, for all $n$ and the norms in $A$ and $\hat{A}$ coincide, it follows that $a_n \xlongrightarrow{n \rightarrow \infty} 0$, and consequently $g(a_n) \longrightarrow g(0) = 0$. Thus, there are no pairs $(m,0)$, with $m \neq 0$ in $\overline{ran\begin{bsmallmatrix}
		-g \\
		i\\
\end{bsmallmatrix}}$. Hence $p((m,0)) = 0 \Longleftrightarrow m = 0$ from which we conclude that $j$ is a monomorphism.\\
$\square$.
\par We present the basic definitions and fundamental facts concerning divisible groups and show one more interesting example of universal monomorphism. Throughout this section, only abelian groups are considered. The category of divisible groups will be denoted by $\textbf{Div}$.
\par An abelian group $(G,+)$ is \textit{divisible} if, for every positive integer $n$ and every $g \in G$, there exists $y \in G$ such that $ny=x$.
\par \textbf{Lemma 3.1.}\cite[Theorem 21.3]{J} Every abelian group $A$ decomposes into a direct sum $A= D \oplus R$ where $D$ is a divisible group and $R$ is a reduced group, i.e., a group containing no nonzero divisible subgroups.
\par \textbf{Lemma 3.2.} A homomorphic image of a divisible group is a divisible group.
\par \textbf{Lemma 3.3.}\cite[Theorem 21.1]{J} Divisible groups are exactly injective objects in the category of abelian groups.
\par \textbf{Lemma 3.4.(The structure theorem for divisible groups)}\cite[Theorem 23.1]{J} For a divisible group $D$, the following decomposition holds:
$$G \cong \bigoplus_{i \in I}\mathbb{Q} \oplus \bigoplus_{p \in P}\mathbb{Z}_{p^\infty}^{\lambda_p}$$
where $$\mathbb{Z}_{p^\infty} = \varinjlim \faktor{\mathbb{Z}}{p^n\mathbb{Z}} = \faktor{\mathbb{Z}[1/p]}{\mathbb{Z}}$$ is the Prüfer group.
\par For an abelian group $A$ denote by $D(A)$ its divisible part.
\par \textbf{Proposition 3.2.} The morphism $f \colon \mathbb{Q} \to \faktor{\mathbb{Q}}{\mathbb{Z}}$ is a bimorphism. Moreover, $f$ is a universal monomorphism but not a universal epimorphism.\\
\textit{Proof.} The morphism $f$ is an epimorphism since it is surjective. Let $i \colon \mathbb{Z}_{p^\infty} \hookrightarrow \faktor{\mathbb{Q}}{\mathbb{Z}}$ be the canonical embedding, where
\[
\faktor{\mathbb{Q}}{\mathbb{Z}} \cong \bigoplus_{p \in P}\mathbb{Z}_{p^\infty}.
\]
Recall that
\[
\mathbb{Z}_{p^\infty} \cong \faktor{\mathbb{Z}[1/p]}{\mathbb{Z}}.
\]
Consider the pullback of $f$ along $i$ in $\mathbf{Ab}$:
$$\begin{tikzcd}
	{\mathbb{Z}[1/p]} \arrow[r, "f'", two heads] \arrow[d] & \faktor{\mathbb{Z}[1/p]}{\mathbb{Z}} \arrow[d, "i", hook] \\
	\mathbb{Q} \arrow[r, "f"', two heads]                  & \faktor{\mathbb{Q}}{\mathbb{Z}}                     
\end{tikzcd}$$
\par Since $\mathbf{Div}$ is a torsion class in the abelian category $\mathbf{Ab}$, taking the divisible part commutes with limits. Indeed, the functor
\[
D \colon \mathbf{Ab} \to \mathbf{Div}, \qquad A \mapsto D(A),
\]
is right adjoint to the canonical inclusion $\mathbf{Div} \hookrightarrow \mathbf{Ab}$.
\par Hence the pullback in $\mathbf{Div}$ is given by $D(\mathbb{Z}[1/p])$. Since $\mathbb{Z}[1/p]$ is not divisible, we have\\ $D(\mathbb{Z}[1/p])=0.$
Consequently, the pullback of $f$ along $i$ in $\mathbf{Div}$ is
\[
0 \to \faktor{\mathbb{Z}[1/p]}{\mathbb{Z}},
\]
which is not an epimorphism. Therefore, $f$ is not a universal epimorphism.
\par On the other hand, clearly,
\[
\Ker_{\mathbf{Ab}} f=\mathbb{Z}.
\]
Since $D(\mathbb{Z})=0$, it follows that
\[
D(\Ker_{\mathbf{Ab}} f)=\Ker_{\mathbf{Div}} f=0,
\]
and hence $f$ is a monomorphism in $\mathbf{Div}$.

Now consider the pushout of $f$ along an arbitrary morphism $g$ in $\mathbf{Ab}$:
$$\begin{tikzcd}
	\mathbb{Z} \arrow[d, "\psi"', two heads, dotted] \arrow[r, hook] & \mathbb{Q} \arrow[d, "g"] \arrow[r, "f"] \arrow[r, "\text{PO}", phantom, bend right=60] & \faktor{\mathbb{Q}}{\mathbb{Z}} \arrow[d] \\
	\Ker_{\mathbf{Ab}} f' \arrow[r, hook]                                      & X \arrow[r, "f'"]                                                                       & Y          
\end{tikzcd}$$
By \cite[Lemma 3]{G}, the morphism $\psi$ is an epimorphism. By the first isomorphism theorem,
\[
\Ker_{\mathbf{Ab}} f' \cong \mathbb{Z}
\qquad\text{or}\qquad
\Ker_{\mathbf{Ab}} f' \cong \faktor{\mathbb{Z}}{n\mathbb{Z}}.
\]
 In either case,
\[
D(\Ker_{\mathbf{Ab}} f')=0,
\]
since both $\mathbb{Z}$ and $\faktor{\mathbb{Z}}{n\mathbb{Z}}$ have trivial divisible part. Therefore,
\[
\Ker_{\mathbf{Div}} f'=0,
\]
so $f'$ is a monomorphism in $\mathbf{Div}$. Thus $f$ is a universal monomorphism.\\
$\square$.
\section{The Main Theorem}
\textbf{Theorem 4.1.} In any quasi-abelian category, Axioms 1 and 1$^\ast$ hold.
\par By the corollary in \cite[p. 193]{R} , every quasi-abelian category is equivalent to a torsion class $\mathcal{T}$ in an abelian category $\mathcal{A}$. Since equivalences preserve limits, colimits,
monomorphisms, and epimorphisms, it is enough to work in such a torsion class.
We now prove an auxiliary statement.
\par \textbf{Proposition 4.1.}
Let $f\colon A\to B$ be a morphism in $\mathcal A$. Every epimorphism
\[
\alpha\colon \Ker f \twoheadrightarrow G
\]
can be realized as the morphism between kernels induced by a suitable pushout.\\
\textit{Proof.}  Let $\alpha : \Ker f \twoheadrightarrow G$ be an epimorphism, and set $N:= \Ker \alpha$. Then there is a monomorphism $i : N \rightarrowtail A$.\\
Consider the pushout of $f$ along $q:=\Coker(i)$. Since the diagram is finite, we may work in a category of modules.
$$\begin{tikzcd}
	N \arrow[d, hook]                                     &                                                                                            &             \\
	\Ker f \arrow[d, "\alpha"', two heads] \arrow[r, hook] & A \arrow[r, "f"] \arrow[d, "q"', two heads] \arrow[r, "\text{PO}", phantom, bend right=67] & B \arrow[d] \\
	G                                                     & \faktor{A}{N} \arrow[r, "f'"]                                                                        & Y          
\end{tikzcd}$$
We claim that
$$q(a) \in \Ker f' \Longleftrightarrow a \in \Ker f$$
For the implication $(\Leftarrow)$, let $a \in \Ker f$. Then $$(f(a),-q(a))=(0,-q(a))$$ hence $-q(a) \in \Ker f'$ and therefore $q(a) \in \Ker f'$.\\
For the converse, assume that $q(a) \in \Ker f'$. Then, by the definition of the pushout, there exists
$b\in A$ such that $$(0,q(a))=(f(b),-q(b)).$$ Hence $b\in \Ker f$ and $q(a+b)=0$, so $a+b\in N$. Since $N\subseteq \Ker f$, it follows that
$a\in \Ker f$.\\
Therefore,
$$\Ker f' = \{q(a) \mid a \in \Ker f\} = q (\Ker f) = \faktor{\Ker f}{N} = G
$$
$\square$.\\
\par Turn to proving Theorem 4.1.
\par We now describe universal monomorphisms in $\mathcal T$. Let $f\colon A\to B$ be a monomorphism in $\mathcal T$. Note that $f$ need not be a monomorphism in $\mathcal A$.
For an arbitrary morphism $g\colon A\to X$ in $\mathcal T$, consider the pushout in $\mathcal A$:
$$\begin{tikzcd}
	\Ker_{\mathcal{A}} f \arrow[r, hook] \arrow[d, "\psi", dashed] & A \arrow[r, "f"] \arrow[d, "g"'] \arrow[r, "\text{PO}", phantom, bend right=60] & B \arrow[d] \\
	\Ker_{\mathcal{A}} f' \arrow[r, hook]                          & X \arrow[r, "f'"']                                                              & Y          
\end{tikzcd} $$
where $\Ker_{\mathcal{A}} f$ is kernel of $f$ in $\mathcal{A}$. Since this is a pushout in an abelian category, $\psi$ is an epimorphism by \cite[Lemma 3]{G}. The monomorphism $f$ is universal if and only if, for every $g$, the induced morphism $f'$ is again
a monomorphism in $\mathcal T$. Equivalently,
\[
F(\Ker_{\mathcal A} f')=0
\]
for every such pushout square, or, in other words,
\[
F\psi=\operatorname{id}_0 .
\] By Proposition 4.1, this yields the following characterization:
\[
f \text{ is universal }
\Longleftrightarrow Fp=\operatorname{id}_0
\text{for every epimorphism } p\colon \Ker_{\mathcal A} f \twoheadrightarrow G.
\]
We are now in a position to verify stability under pullbacks. Let $\delta\colon Z\to B$ be an arbitrary morphism in $\mathcal T$, and consider the pullback of $f$ along $\delta$ in
$\mathcal A$:
$$\begin{tikzcd}
	\Ker_{\mathcal{A}}Ff' \arrow[r, hook] \arrow[d, "r", tail]          & F(D) \arrow[d, "\mu", tail] \arrow[rd, "Ff'"]                                             &                       \\
	\Ker_{\mathcal{A}}f' \arrow[r, hook] \arrow[d, equal] & D \arrow[r, "f'"] \arrow[d] \arrow[r, "\text{PB}", phantom, bend right=49, shift right=2] & Z \arrow[d, "\delta"] \\
	\Ker_{\mathcal{A}}f \arrow[r, hook]                                 & A \arrow[r, "f"']                                                                         & B                    
\end{tikzcd}$$
As observed above, the pullback in $\mathcal T$ is obtained by applying $F$.\\
Since
\[
f'\circ \mu \circ \ker_{\mathcal A}Ff' = Ff'\circ \ker_{\mathcal A}Ff' = 0,
\]
there exists a unique morphism
\[
r\colon \Ker_{\mathcal A}(Ff')\rightarrowtail \Ker_{\mathcal A} f'
\]
such that
\[
(\ker_{\mathcal A} f')\circ r=\mu\circ \ker_{\mathcal A}Ff' .
\]
Since $\mu\circ \ker_{\mathcal A}Ff'$ is a monomorphism, so is $r$. To prove that $Ff'$ is universal, consider an arbitrary pushout of $Ff'$ along a morphism $h\colon \Ker_{\mathcal A}(Ff')\to R$. Computing this pushout successively --- first along $\mu$ and then along $f'$ --- we obtain a diagram
in which the right-hand pushout is associated with $f'$:
$$\begin{tikzcd}
	\Ker_{\mathcal{A}}Ff' \arrow[r, hook] \arrow[d, "\psi'"', dashed] & F(D) \arrow[r, "\mu", tail] \arrow[d, "h"] \arrow[rr, "Ff'", bend left] \arrow[r, "\text{PO}", phantom, bend right=60] & D \arrow[r, "f'"] \arrow[d] \arrow[r, "\text{PO}", phantom, bend right=60] & Z \arrow[d] \\
	\Ker_{A}\xi \arrow[r, hook]                                       & R \arrow[r, "\mu'"', tail] \arrow[rr, "\xi"', bend right]                                                              & R' \arrow[r, "\theta"']                                                    & W          
\end{tikzcd}$$
Consider the right-hand pushout:
$$\begin{tikzcd}
	\Ker_{\mathcal{A}}f \arrow[r, hook] \arrow[d, "\varphi"', dashed] & D \arrow[r, "f'"] \arrow[d] \arrow[r, "\text{PO}", phantom, bend right=60] & Z \arrow[d] \\
	\Ker_{\mathcal{A}}\theta \arrow[r, hook]                       & R' \arrow[r, "\theta"']                                                    & W          
\end{tikzcd}$$

Since $f$ is universal, the characterization above implies that $F\varphi=\operatorname{id}_0,$ hence $F(\Ker_{\mathcal A}\theta)=0.$ Arguing as above, we get a monomorphism
\[
\zeta\colon \Ker_{\mathcal A}\xi \rightarrowtail \Ker_{\mathcal A}\theta .
\]
Since monomorphisms in $\mathcal A$ are precisely kernels and $F$ preserves limits, it follows that
$F\zeta$ is again a kernel. But $F(\Ker_{\mathcal A}\theta)=0$, and therefore
\[
F\zeta=\operatorname{id}_0.
\]
Consequently,
\[
F(\Ker_{\mathcal A}\xi)=0,
\]
which shows that $Ff'$ is universal mono.
Thus Axiom 1 holds. Axiom 1$^\ast$ follows by duality. This completes the proof of the theorem.

\section{Universal mono and epimorpisms in category of Banach spaces.}
First, we demonstrate the equality $M_u=M_c$. The fact that the category $\textbf{Ban}$ is quasi-abelian and non-integral \cite[Theorem 3.2]{A}.
\par \textbf{Lemma 5.1.}\cite[Proposition 3.2.16]{B} In the category of normed and Banach spaces any morphism has a decomposition of this form:
$$\begin{tikzcd}
	u^{-1}(0) \arrow[r, "\ker u"] & E \arrow[r, "u"] \arrow[d, "\coim u"'] & F \arrow[r, "\coker u"]             & \faktor{F}{\overline{u(E)}} \\
	& \faktor{E}{u^{-1}(0)} \arrow[r, "\overline u"']  & \overline{u(E)} \arrow[u, "\im u"'] &                  
\end{tikzcd}$$
A continuous linear operator $u: X \to Y$ between normed spaces is said to be \textit{bounded below} if there exists a positive $\delta$ such that $\|u(x)\| \geq \delta \|x\|$ for all $x \in X$. Note that such an operator is injective: if $u(a)=u(b)$, then $u(a-b)=0$, and boundedness below yields $0 = \|u(a-b)\| \geq \delta \|a-b\|$, whence $a-b = 0$ and $a = b$. Thus a continuous linear operator that is bounded below is a monomorphism.

\par \textbf{Lemma 5.2.} A monomorphism is bounded below if and only if it is open onto its range.\\
\textit{Proof.} Let $u: X \to Y$ be a monomorphism, that is, an injective continuous linear operator. Then $\|u(x)\| \geq \delta \|x\|$ if and only if $\|x\| \leq \frac{1}{\delta}\|u(x)\|$, i.e. $u^{-1}: u(X) \to X$ is continuous, which is equivalent to $u: X \to u(X)$ being a homeomorphism and hence an open mapping.\\
$\square$.
\par Observe that in the category of Banach spaces, a continuous linear operator is bounded below if and only if it is injective and has closed range. The fact that the range of an operator bounded below is closed is proved in \cite[p. 21]{M}. We prove the converse: if an operator is injective and its range is closed, then the operator is bounded below. Indeed, let $u: X \to Y$ be an injective continuous linear operator between Banach spaces with $u(X)$ closed; then $u(X)$ is a Banach subspace of $Y$ and $u: X \to u(X)$ is bijective. By a corollary of the Banach Open Mapping Theorem, $u^{-1}: u(X) \to X$ is a continuous linear operator, so $\|u^{-1}(u(x))\| \leq C\|u(x)\|$, which implies $\|u(x)\| \geq \frac{\|x\|}{C}$. Since in $\textbf{Ban}$ the strict morphisms are precisely those with closed range \cite[Proposition 3.2.16 (e)]{B}, this will imply the desired equality of classes.
\par Let us recall some fundamental theorems from functional analysis:
\par \textbf{Lemma 5.3.}\cite[p. 21]{M} $u: X \to Y$ is not bounded below if and only if there is a sequence of unit vectors $\{x_n\}$ in $X$, such that $u(x_n)\xrightarrow{n \to \infty} 0$.
\par Let $u : X \to Y$ be a monomorphism that is not bounded below. Then by Lemma 5.3, there exists a sequence $\{x_n\}$ in $X$ with $\|x_n\|_X=1$ such that $u(x_n)\xrightarrow{n \to \infty} 0$. Define $p(x):=\|u(x)\|_Y$ and for those $n$, where $p(x_n) \neq 0$, set $y_n:= \frac{x_n}{p(x_n)}$. Then $p(y_n)=1$ and $\|y_n\|_X = \frac{1}{p(x_n)} \xrightarrow{n \to \infty} \infty$. Hence the set $S:= \{ x \in X : p(x) \leq 1 \}$ is unbounded in $X$.
\par Let $X^*$ be the space of continuous linear functionals $f: X \to \mathbb{K}$ here $\mathbb{K}$ is the base field. Fix $x \in S$ and consider the operator $T_x : X^* \to \mathbb{K}$ defined by the rule $T_x(\psi) = \psi(x)$. It is well-known that $T_x$ is a continuous linear functional on $X^*$ with $\|T_x\|=\|x\|$.
\par \textbf{Lemma 5.4.} There exists $\varphi \in X^*$ such that $\sup \{ |\varphi(x)| : p(x)\leq 1 \} = \infty$.\\
\textit{Proof.} Suppose on the contrary that for every $\psi \in X^*$, $\sup \{ |\psi(x)| : x \in S \} < \infty$. Then $\sup \{ |T_x (\psi)| : x \in S \} < \infty$; by the Uniform Bounded Principle we obtain $\sup \{ \|T_x\| : x \in S \} < \infty$. Since $\|T_x\|=\|x\|$, we conclude that $\sup \{ \|x\| : x \in S \} < \infty$ which is a contradiction since $S$ is unbounded.\\
$\square$.
\par From the condition $\sup \{ |\varphi(x)| : p(x)\leq 1 \} = \infty$ we obtain that for every $n \in \mathbb{N}$ there exists $a_n \in X$ such that $p(a_n) \leq 1$ and $|\varphi(a_n)| \geq n$. Define $b_n:= \frac{a_n}{\varphi(a_n)}$; then $\varphi(b_n)=1$ and $p(b_n) \leq \frac{1}{n}$.\\
Consider the pushout of $u$ along $\varphi$:

	$$\begin{tikzcd}
		X \arrow[rr, "u"] \arrow[dd, "\varphi"'] \arrow[rr, "\text{PO}", phantom, bend right=60] &                                    & Y \arrow[d, "{(0,id)}"]            \\
		&                                    & \mathbb{K} \oplus Y \arrow[d, "p"] \\
		\mathbb{K} \arrow[r, "{(id,0)}"]                                                         & \mathbb{K} \oplus Y \arrow[r, "p"] & D                                 
	\end{tikzcd}$$
Here $D := \faktor{\mathbb{K} \oplus Y}{\overline{ran\begin{bsmallmatrix}
			-u \\
			\varphi\\
\end{bsmallmatrix}}}$ $p$ is the canonical projection.
\par Let us denote $j :=p \circ (id,0)$. By definition, $j(m) = p((m,0))$, i.e., $j(m)=0$ if $(m,0) \in \overline{ran\begin{bsmallmatrix}
		-u \\
		\varphi\\
\end{bsmallmatrix}}$.
Show that $(1,0) \in \overline{ran\begin{bsmallmatrix}
		-u \\
		\varphi\\
\end{bsmallmatrix}}$.
$\begin{bsmallmatrix}
	-u \\
	\varphi\\
\end{bsmallmatrix} (b_n) = \begin{bsmallmatrix}
	-u(b_n) \\
	\varphi(b_n) \\
\end{bsmallmatrix} \xrightarrow{n \to \infty} \begin{bsmallmatrix}
	0 \\
	1\\
\end{bsmallmatrix}$. Hence $j(1)=0$, from which we conclude that 
$j$ is not a monomorphism.
\par Thus, we have demonstrated that universal monomorphisms in $\textbf{Ban}$ are exactly the operators bounded below, and since these are strict, we obtain $M_u=M_c$. Now prove that $P_u=P_c$.
\par \textbf{Remark:} The following lemma is a folklore but we have not found this fact in the literature.
\par \textbf{Lemma 5.5.} Strict epimorphisms in $\textbf{Ban}$ are preciesly surjective operators.\\
\textit{Proof.} Let $f: X \to Y$ be a surjection. Then $f(X)=Y$, hence its image is both dense and closed, and consequently the morphism is both an epimorphism and strict.
\par For the converse, assume that $f$ is a strict epimorphism. Then $\Coker f = 0$. On the other hand, by Lemma 5.1,
\[
\Coker f = \faktor{B}{\overline{f(A)}}.
\]
Hence $\faktor{B}{\overline{f(A)}} = 0.$
Since $f$ is strict, $f(A)$ is closed, and so $\overline{f(A)} = f(A)$. Therefore, $\faktor{B}{f(A)} = 0$, and thus $B = f(A)$.
$\square$.
\par \textbf{Proposition 5.1.} $f: X \to Y$ is a universal epimorphism in \textbf{Ban} if and only if $f$ is a surjection.\\
\textit{Proof.} $(\Longleftarrow)$ Obviously by the axioms of quasi-abelianity.\\
$(\Longrightarrow)$ Suppose $f: X \to Y$ is not surjective. Choose an element $y_0 \in Y\setminus f(X)$. Define $\varphi : \mathbb{K} \to Y$ by the rule $\varphi(\lambda):=\lambda y_0$.
Consider the following pullback:\\
$$\begin{tikzcd}
	D \arrow[r, "\pi"] \arrow[d] \arrow[r, "\text{PB}", phantom, bend right=60] & \mathbb{K} \arrow[d, "\varphi"] \\
	X \arrow[r, "f"']                                                           & Y                              
\end{tikzcd}$$
Here $D=\{(a,\lambda):f(a)=\varphi(\lambda)\}$. However, $f(a)=\lambda y_0$ is possible only when $\lambda = 0$. Then $\pi(D)=\{0\}$ is not dense in $\mathbb{K}$ and hence is not an epimorphism.\\
$\square$.
\par Since in $\textbf{Ban}$ we have $M_u=M_c$ and $P_u=P_c$, it is evident by the axioms of quasi-abelianity that Axioms 1 and 1* also hold. Note that in the proof of $P_u=P_c$, we did not use completeness anywhere; hence, this equality of classes also holds in $\textbf{Norm}$. For Banach spaces, we have essentially shown that universal monomorphisms are precisely the operators bounded below. However, $X^*$ is a Banach space for any normed space $X$ and completeness was not required elsewhere except in Lemma 5.3, which easily generalizes to normed spaces. Thus, we obtain that universal monomorphisms in $\textbf{Norm}$ are exactly the operators bounded below (that is, the injective operators open onto their range).
\section*{Acknowledgments}
The author would like to express sincere gratitude to Professor Yaroslav Kopylov for his guidance, helpful discussions, and continued support.

\end{document}